\documentclass[12pt]{article}
\usepackage{amsmath,amssymb,amsthm,latexsym,amscd}

\title{Combinations related to classes \\ of finite and countably categorical structures \\ and their theories\footnote{
\newline\indent \ \ \ The research is partially supported by the Grants Council
(under RF President) for State Aid of Leading Scientific Schools
(Grant NSh-6848.2016.1) and by Committee of Science in Education
and Science Ministry of the Republic of Kazakhstan (Grant No.
0830/GF4). } }
\author{Sergey V.
Sudoplatov\footnote{sudoplat@math.nsc.ru}}
\date{}
\begin{document}
\maketitle

\begin{abstract}
We consider and characterize classes of finite and countably
categorical structures and their theories preserved under
$E$-operators and $P$-operators. We describe $e$-spectra and
families of finite cardinalities for structures belonging to
closures with respect to $E$-operators and $P$-operators.

\end{abstract}

We continue to study structural properties of $E$-combin\-a\-tions
and $P$-combin\-a\-tions of structures and their theories
\cite{cs, cl, lut, rest, lft} applying the general context to the
classes of finitely categorical and $\omega$-categorical theories.

Approximations of structures by finite ones as well as
correspondent approximations of theories were studied in a series
of papers, e.g. \cite{Rosen, Vaa, Zil10}. We consider these
approximations in the context of structural combinations.

We consider and describe $e$-spectra and families of finite
cardinalities for structures belonging to closures with respect to
$E$-operators and $P$-operators.

\section{Preliminaries}

Throughout the paper we use the following terminology in \cite{cs,
cl} as well as in \cite{SuGP, SuMCT}.

Let $P=(P_i)_{i\in I}$, be a family of nonempty unary predicates,
$(\mathcal{A}_i)_{i\in I}$ be a family of structures such that
$P_i$ is the universe of $\mathcal{A}_i$, $i\in I$, and the
symbols $P_i$ are disjoint with languages for the structures
$\mathcal{A}_j$, $j\in I$. The structure
$\mathcal{A}_P\rightleftharpoons\bigcup\limits_{i\in
I}\mathcal{A}_i$\index{$\mathcal{A}_P$} expanded by the predicates
$P_i$ is the {\em $P$-union}\index{$P$-union} of the structures
$\mathcal{A}_i$, and the operator mapping $(\mathcal{A}_i)_{i\in
I}$ to $\mathcal{A}_P$ is the {\em
$P$-operator}\index{$P$-operator}. The structure $\mathcal{A}_P$
is called the {\em $P$-combination}\index{$P$-combination} of the
structures $\mathcal{A}_i$ and denoted by ${\rm
Comb}_P(\mathcal{A}_i)_{i\in I}$\index{${\rm
Comb}_P(\mathcal{A}_i)_{i\in I}$} if
$\mathcal{A}_i=(\mathcal{A}_P\upharpoonright
A_i)\upharpoonright\Sigma(\mathcal{A}_i)$, $i\in I$. Structures
$\mathcal{A}'$, which are elementary equivalent to ${\rm
Comb}_P(\mathcal{A}_i)_{i\in I}$, will be also considered as
$P$-combinations.

Clearly, all structures $\mathcal{A}'\equiv {\rm
Comb}_P(\mathcal{A}_i)_{i\in I}$ are represented as unions of
their restrictions $\mathcal{A}'_i=(\mathcal{A}'\upharpoonright
P_i)\upharpoonright\Sigma(\mathcal{A}_i)$ if and only if the set
$p_\infty(x)=\{\neg P_i(x)\mid i\in I\}$ is inconsistent. If
$\mathcal{A}'\ne{\rm Comb}_P(\mathcal{A}'_i)_{i\in I}$, we write
$\mathcal{A}'={\rm Comb}_P(\mathcal{A}'_i)_{i\in
I\cup\{\infty\}}$, where
$\mathcal{A}'_\infty=\mathcal{A}'\upharpoonright
\bigcap\limits_{i\in I}\overline{P_i}$, maybe applying
Morleyzation. Moreover, we write ${\rm
Comb}_P(\mathcal{A}_i)_{i\in I\cup\{\infty\}}$\index{${\rm
Comb}_P(\mathcal{A}_i)_{i\in I\cup\{\infty\}}$} for ${\rm
Comb}_P(\mathcal{A}_i)_{i\in I}$ with the empty structure
$\mathcal{A}_\infty$.

Note that if all predicates $P_i$ are disjoint, a structure
$\mathcal{A}_P$ is a $P$-combination and a disjoint union of
structures $\mathcal{A}_i$. In this case the $P$-combination
$\mathcal{A}_P$ is called {\em
disjoint}.\index{$P$-combination!disjoint} Clearly, for any
disjoint $P$-combination $\mathcal{A}_P$, ${\rm
Th}(\mathcal{A}_P)={\rm Th}(\mathcal{A}'_P)$, where
$\mathcal{A}'_P$ is obtained from $\mathcal{A}_P$ replacing
$\mathcal{A}_i$ by pairwise disjoint
$\mathcal{A}'_i\equiv\mathcal{A}_i$, $i\in I$. Thus, in this case,
similar to structures the $P$-operator works for the theories
$T_i={\rm Th}(\mathcal{A}_i)$ producing the theory $T_P={\rm
Th}(\mathcal{A}_P)$\index{$T_P$}, being {\em
$P$-combination}\index{$P$-combination} of $T_i$, which is denoted
by ${\rm Comb}_P(T_i)_{i\in I}$.\index{${\rm Comb}_P(T_i)_{i\in
I}$} In general, for non-disjoint case, the theory $T_P$ will be
also called a {\em $P$-combination}\index{$P$-combination} of the
theories $T_i$, but in such a case we will keep in mind that this
$P$-combination is constructed with respect (and depending) to the
structure $\mathcal{A}_P$, or, equivalently, with respect to
any/some $\mathcal{A}'\equiv\mathcal{A}_P$.

For an equivalence relation $E$ replacing disjoint predicates
$P_i$ by $E$-classes we get the structure
$\mathcal{A}_E$\index{$\mathcal{A}_E$} being the {\em
$E$-union}\index{$E$-union} of the structures $\mathcal{A}_i$. In
this case the operator mapping $(\mathcal{A}_i)_{i\in I}$ to
$\mathcal{A}_E$ is the {\em $E$-operator}\index{$E$-operator}. The
structure $\mathcal{A}_E$ is also called the {\em
$E$-combination}\index{$E$-combination} of the structures
$\mathcal{A}_i$ and denoted by ${\rm Comb}_E(\mathcal{A}_i)_{i\in
I}$\index{${\rm Comb}_E(\mathcal{A}_i)_{i\in I}$}; here
$\mathcal{A}_i=(\mathcal{A}_E\upharpoonright
A_i)\upharpoonright\Sigma(\mathcal{A}_i)$, $i\in I$. Similar
above, structures $\mathcal{A}'$, which are elementary equivalent
to $\mathcal{A}_E$, are denoted by ${\rm
Comb}_E(\mathcal{A}'_j)_{j\in J}$, where $\mathcal{A}'_j$ are
restrictions of $\mathcal{A}'$ to its $E$-classes. The
$E$-operator works for the theories $T_i={\rm Th}(\mathcal{A}_i)$
producing the theory $T_E={\rm Th}(\mathcal{A}_E)$\index{$T_E$},
being {\em $E$-combination}\index{$E$-combination} of $T_i$, which
is denoted by ${\rm Comb}_E(T_i)_{i\in I}$\index{${\rm
Comb}_E(T_i)_{i\in I}$} or by ${\rm
Comb}_E(\mathcal{T})$\index{${\rm Comb}_E(\mathcal{T})$}, where
$\mathcal{T}=\{T_i\mid i\in I\}$.

Clearly, $\mathcal{A}'\equiv\mathcal{A}_P$ realizing $p_\infty(x)$
is not elementary embeddable into $\mathcal{A}_P$ and can not be
represented as a disjoint $P$-combination of
$\mathcal{A}'_i\equiv\mathcal{A}_i$, $i\in I$. At the same time,
there are $E$-combinations such that all
$\mathcal{A}'\equiv\mathcal{A}_E$ can be represented as
$E$-combinations of some $\mathcal{A}'_j\equiv\mathcal{A}_i$. We
call this representability of $\mathcal{A}'$ to be the {\em
$E$-representability}.

If there is $\mathcal{A}'\equiv\mathcal{A}_E$ which is not
$E$-representable, we have the $E'$-representability replacing $E$
by $E'$ such that $E'$ is obtained from $E$ adding equivalence
classes with models for all theories $T$, where $T$ is a theory of
a restriction $\mathcal{B}$ of a structure
$\mathcal{A}'\equiv\mathcal{A}_E$ to some $E$-class and
$\mathcal{B}$ is not elementary equivalent to the structures
$\mathcal{A}_i$. The resulting structure $\mathcal{A}_{E'}$ (with
the $E'$-representability) is a {\em
$e$-completion}\index{$e$-completion}, or a {\em
$e$-saturation}\index{$e$-saturation}, of $\mathcal{A}_{E}$. The
structure $\mathcal{A}_{E'}$ itself is called {\em
$e$-complete}\index{Structure!$e$-complete}, or {\em
$e$-saturated}\index{Structure!$e$-saturated}, or {\em
$e$-universal}\index{Structure!$e$-universal}, or {\em
$e$-largest}\index{Structure!$e$-largest}.

For a structure $\mathcal{A}_E$ the number of {\em
new}\index{Structure!new} structures with respect to the
structures $\mathcal{A}_i$, i.~e., of the structures $\mathcal{B}$
which are pairwise elementary non-equivalent and elementary
non-equivalent to the structures $\mathcal{A}_i$, is called the
{\em $e$-spectrum}\index{$e$-spectrum} of $\mathcal{A}_E$ and
denoted by $e$-${\rm Sp}(\mathcal{A}_E)$.\index{$e$-${\rm
Sp}(\mathcal{A}_E)$} The value ${\rm sup}\{e$-${\rm
Sp}(\mathcal{A}'))\mid\mathcal{A}'\equiv\mathcal{A}_E\}$ is called
the {\em $e$-spectrum}\index{$e$-spectrum} of the theory ${\rm
Th}(\mathcal{A}_E)$ and denoted by $e$-${\rm Sp}({\rm
Th}(\mathcal{A}_E))$.\index{$e$-${\rm Sp}({\rm
Th}(\mathcal{A}_E))$}

If $\mathcal{A}_E$ does not have $E$-classes $\mathcal{A}_i$,
which can be removed, with all $E$-classes
$\mathcal{A}_j\equiv\mathcal{A}_i$, preserving the theory ${\rm
Th}(\mathcal{A}_E)$, then $\mathcal{A}_E$ is called {\em
$e$-prime}\index{Structure!$e$-prime}, or {\em
$e$-minimal}\index{Structure!$e$-minimal}.

For a structure $\mathcal{A}'\equiv\mathcal{A}_E$ we denote by
${\rm TH}(\mathcal{A}')$ the set of all theories ${\rm
Th}(\mathcal{A}_i)$\index{${\rm Th}(\mathcal{A}_i)$} of
$E$-classes $\mathcal{A}_i$ in $\mathcal{A}'$.

By the definition, an $e$-minimal structure $\mathcal{A}'$
consists of $E$-classes with a minimal set ${\rm
TH}(\mathcal{A}')$. If ${\rm TH}(\mathcal{A}')$ is the least for
models of ${\rm Th}(\mathcal{A}')$ then $\mathcal{A}'$ is called
{\em $e$-least}.\index{Structure!$e$-least}

\medskip
{\bf Definition} \cite{cl}. Let $\overline{\mathcal{T}}$ be the
class of all complete elementary theories of relational languages.
For a set $\mathcal{T}\subset\overline{\mathcal{T}}$ we denote by
${\rm Cl}_E(\mathcal{T})$ the set of all theories ${\rm
Th}(\mathcal{A})$, where $\mathcal{A}$ is a structure of some
$E$-class in $\mathcal{A}'\equiv\mathcal{A}_E$,
$\mathcal{A}_E={\rm Comb}_E(\mathcal{A}_i)_{i\in I}$, ${\rm
Th}(\mathcal{A}_i)\in\mathcal{T}$. As usual, if $\mathcal{T}={\rm
Cl}_E(\mathcal{T})$ then $\mathcal{T}$ is said to be {\em
$E$-closed}.\index{Set!$E$-closed}

The operator ${\rm Cl}_E$ of $E$-closure can be naturally extended
to the classes $\mathcal{T}\subset\overline{\mathcal{T}}$ as
follows: ${\rm Cl}_E(\mathcal{T})$ is the union of all ${\rm
Cl}_E(\mathcal{T}_0)$ for subsets
$\mathcal{T}_0\subseteq\mathcal{T}$.

For a set $\mathcal{T}\subset\overline{\mathcal{T}}$ of theories
in a language $\Sigma$ and for a sentence $\varphi$ with
$\Sigma(\varphi)\subseteq\Sigma$ we denote by
$\mathcal{T}_\varphi$\index{$\mathcal{T}_\varphi$} the set
$\{T\in\mathcal{T}\mid\varphi\in T\}$.

\medskip
{\bf Proposition 1.1} \cite{cl}. {\em If
$\mathcal{T}\subset\overline{\mathcal{T}}$ is an infinite set and
$T\in\overline{\mathcal{T}}\setminus\mathcal{T}$ then $T\in{\rm
Cl}_E(\mathcal{T})$ {\rm (}i.e., $T$ is an {\sl accumulation
point} for $\mathcal{T}$ with respect to $E$-closure ${\rm
Cl}_E${\rm )} if and only if for any formula $\varphi\in T$ the
set $\mathcal{T}_\varphi$ is infinite.}

\medskip
{\bf Theorem 1.2} \cite{cl}. {\em If $\mathcal{T}'_0$ is a
generating set for a $E$-closed set $\mathcal{T}_0$ then the
following conditions are equivalent:

$(1)$ $\mathcal{T}'_0$ is the least generating set for
$\mathcal{T}_0$;

$(2)$ $\mathcal{T}'_0$ is a minimal generating set for
$\mathcal{T}_0$;

$(3)$ any theory in $\mathcal{T}'_0$ is isolated by some set
$(\mathcal{T}'_0)_\varphi$, i.e., for any $T\in\mathcal{T}'_0$
there is $\varphi\in T$ such that
$(\mathcal{T}'_0)_\varphi=\{T\}$;

$(4)$ any theory in $\mathcal{T}'_0$ is isolated by some set
$(\mathcal{T}_0)_\varphi$, i.e., for any $T\in\mathcal{T}'_0$
there is $\varphi\in T$ such that
$(\mathcal{T}_0)_\varphi=\{T\}$.}

\medskip
{\bf Definition} \cite{cl}. For a set
$\mathcal{T}\subset\overline{\mathcal{T}}$ we denote by ${\rm
Cl}_P(\mathcal{T})$\index{${\rm Cl}_P(\mathcal{T})$} the set of
all theories ${\rm Th}(\mathcal{A})$ such that ${\rm
Th}(\mathcal{A})\in\mathcal{T}$ or $\mathcal{A}$ is a structure of
type $p_\infty(x)$ in $\mathcal{A}'\equiv\mathcal{A}_P$, where
$\mathcal{A}_P={\rm Comb}_P(\mathcal{A}_i)_{i\in I}$ and ${\rm
Th}(\mathcal{A}_i)\in\mathcal{T}$ are pairwise distinct. As above,
if $\mathcal{T}={\rm Cl}_P(\mathcal{T})$ then $\mathcal{T}$ is
said to be {\em $P$-closed}.\index{Set!$P$-closed}

Using above only disjoint $P$-combinations $\mathcal{A}_P$ we get
the closure ${\rm Cl}^d_P(\mathcal{T})$\index{${\rm
Cl}^d_P(\mathcal{T})$} being a subset of ${\rm
Cl}_P(\mathcal{T})$.

The closure operator ${\rm Cl}^{d,r}_P$\index{${\rm Cl}^{d,r}_P$}
is obtained from ${\rm Cl}^{d}_P$ permitting repetitions of
structures for predicates $P_i$.

\medskip
Replacing $E$-classes by unary predicates $P_i$ (not necessary
disjoint) being universes for structures $\mathcal{A}_i$ and
restricting models of ${\rm Th}(\mathcal{A}_P)$ to the set of
realizations of $p_\infty(x)$ we get the {\em
$e$-spectrum}\index{$e$-spectrum} $e$-${\rm Sp}({\rm
Th}(\mathcal{A}_P))$\index{$e$-${\rm Sp}({\rm
Th}(\mathcal{A}_P))$}, i.~e., the number of pairwise elementary
non-equivalent restrictions of $\mathcal{M}\models{\rm
Th}(\mathcal{A}_P)$ to $p_\infty(x)$ such that these restrictions
are not elementary equivalent to the structures $\mathcal{A}_i$.

\medskip
{\bf Definition} \cite{SuGP, SuMCT}. A {\em $n$-dimensional
cube},\index{Cube!$n$-dimensional} or a {\em
$n$-cube}\index{$n$-cube} (where $n\in\omega$) is a graph
isomorphic to the graph ${\cal Q}_n$\index{${\cal Q}_n$} with the
universe $\{0,1\}^n$ and such that any two vertices
$(\delta_1,\ldots,\delta_n)$ and $(\delta'_1,\ldots,\delta'_n)$
are adjacent if and only if these vertices differ exactly in one
coordinate. The described graph ${\cal Q}_n$ is called the {\em
canonical representative}\index{Representative!canonical} for the
class of $n$-cubes.

Let $\lambda$ be an infinite cardinal. A {\em
$\lambda$-dimensional cube},\index{Cube!$\lambda$-dimensional} or
a {\em $\lambda$-cube},\index{$\lambda$-cube} is a graph
isomorphic to a graph $\Gamma=\langle X;R\rangle$ satisfying the
following conditions:

(1) the universe $X\subseteq\{0,1\}^\lambda$ is generated from an
arbitrary function $f\in X$ by the operator $\langle
f\rangle$\index{$\langle f\rangle$} attaching, to the set $\{f\}$,
all results of substitutions for any finite tuples
$(f(i_1),\ldots,f(i_m))$ by tuples $(1-f(i_1),\ldots,1-f(i_m))$;

(2) the relation $R$ consists of edges connecting functions
differing exactly in one coordinate (the ({\em $i$-th}) {\em
coordinate}\index{Coordinate} of function $g\in\{0,1\}^\lambda$ is
the value $g(i)$ correspondent to the argument $i<\lambda$).

The described graph ${\cal Q}\rightleftharpoons{\cal
Q}_f$\index{${\cal Q}_f$} with the universe $\langle f\rangle$ is
a {\em canonical representative} for the class of $\lambda$-cubes.

Note that the canonical representative of the class of $n$-cubes
(as well as the canonical representatives of the class of
$\lambda$-cubes) are generated by any its function:
$\{0,1\}^n=\langle f\rangle$, where $f\in\{0,1\}^n$. Therefore the
universes of canonical representatives ${\cal Q}_f$ of $n$-cubes
like $\lambda$-cubes, is denoted by $\langle f\rangle$.

\section{Closed classes of finitely categorical \\ and $\omega$-categorical theories}

Remind that a countable complete theory $T$ is {\em
$\omega$-categorical} if $T$ has exactly one countable model up to
isomorphisms, i.e. $I(T,\omega)=1$. A countable theory $T$ is {\em
$n$-categorical}, for natural $n\geq 1$, if $T$ has exactly one
$n$-element model up to isomorphisms, i.e. $I(T,n)=1$. A countable
theory $T$ is {\em finitely categorical} if $T$ is $n$-categorical
for some $n\in\omega\setminus\{0\}$.

The classes of all finitely and $\omega$-categorical theories will
be denoted by $\overline{\mathcal{T}}_{\rm fin}$ and
$\overline{\mathcal{T}}_{\omega,1}$, respectively.

\medskip
Let $\mathcal{T}$ be a set (class) of 
theories in $\overline{\mathcal{T}}_{\rm
fin}\cup\overline{\mathcal{T}}_{\omega,1}$, $T$ be a theory in
$\mathcal{T}$. By Ryll-Nardzewski theorem, $S^n(T)$ is finite for
any $n$. Then, for any $n$, classes
$\lceil\varphi(\bar{x})\rceil=\{\varphi'(\bar{x})\mid\varphi'(\bar{x})\equiv\varphi(\bar{x})\}$
of $T$-formulas with $n$ free variables and
$\lceil\varphi(\bar{x})\rceil\leq\lceil\psi(\bar{x})\rceil\Leftrightarrow\varphi(\bar{x})\vdash\psi(\bar{x})$
form a finite Boolean algebra $\mathcal{B}_n(T)$ with $2^{m_n}$
elements, where $m_n$ is the number of $n$-types of $T$.

The algebra $\mathcal{B}_n(T)$ can be interpreted as a $m_n$-cube
$\mathcal{C}_{m_n}(T)$, whose vertices form the universe $B_n(T)$
of $\mathcal{B}_n(T)$, edges $[a,b]$ link vertices $a$ and $b$
such that $a$ $\leqslant$-covers $b$ or $b$ $\leqslant$-covers
$a$, and each vertex $a$ is marked by some
$u_a\rightleftharpoons\lceil\varphi(\bar{x})\rceil$, where
$a\leqslant b\Leftrightarrow u_a\leq u_b$. The label $0$ is used
for the vertex corresponding to
$\lceil\neg\bar{x}\approx\bar{x}\rceil$ and $1$~--- for the vertex
corresponding to $\lceil\bar{x}\approx\bar{x}\rceil$.

Obviously, the sets $\lceil\varphi(\bar{x})\rceil$ and the
relation $\leq$ depend on the theory $T$ but we omit $T$ if the
theory is fixed or it is clear by the context.

Clearly, algebras $\mathcal{B}_n(T_1)$ and $\mathcal{B}_n(T_2)$,
for $T_1,T_2\in \mathcal{T}$, may be not coordinated: it is
possible $\lceil\varphi(\bar{x})\rceil<\lceil\psi(\bar{x})\rceil$
for $T_1$ whereas
$\lceil\psi(\bar{x})\rceil<\lceil\varphi(\bar{x})\rceil$ for
$T_2$. If $\lceil\varphi(\bar{x})\rceil<\lceil\psi(\bar{x})\rceil$
for $T_1$ and
$\lceil\varphi(\bar{x})\rceil<\lceil\psi(\bar{x})\rceil$ for
$T_2$, we say that $T_2$ {\em witnesses} that
$\lceil\varphi(\bar{x})\rceil<\lceil\psi(\bar{x})\rceil$ for $T_1$
(and vice versa).

At the same time, if a countable theory $T_0$ does not belong to
$\overline{\mathcal{T}}_{\rm
fin}\cup\overline{\mathcal{T}}_{\omega,1}$ then for some $n\geq
1$, $\mathcal{B}_n(T_0)$ is infinite and therefore there is a
formula $\varphi(\bar{x})$, for instance $(\bar x\approx\bar x)$,
such that for the label $u=\lceil\varphi(\bar{x})\rceil$ there is
an infinite decreasing chain $(u_k)_{k\in\omega}$ of labels:
$u_{k+1}<u_k<u$, witnessed by some formulas $\varphi_k(\bar{x})$.
In such a case, if $T_0\in{\rm Cl}_E(\mathcal{T})$, then by
Proposition 1.1 for any finite sequence $(u_l,\ldots,u_0,u)$ there
are infinitely many theories in $\mathcal{T}$ witnessing that
$u_l<\ldots < u_0<u$. In particular, cardinalities $m_n$ for
Boolean algebras $\mathcal{B}_n(T)$ and for cubes
$\mathcal{C}_{m_n}(T)$ are unbounded for $\mathcal{T}$: distances
$\rho_{n,T}(0,u)$ are unbounded for the cubes
$\mathcal{C}_{m_n}(T)$, i.e., ${\rm sup}\{\rho_{n,T}(0,u)\mid
T\in\mathcal{T}\}=\infty$. It is equivalent to take
$(\bar{x}\approx\bar{x}$) for $\varphi(\bar{x})$ and to get ${\rm
sup}\{\rho_{n,T}(0,1)\mid T\in\mathcal{T}\}=\infty$.

Thus we get the following

\medskip
{\bf Theorem 2.1.} {\em Let $\mathcal{T}$ be a class of theories
in $\overline{\mathcal{T}}_{\rm
fin}\cup\overline{\mathcal{T}}_{\omega,1}$. The following
conditions are equivalent:

$(1)$ ${\rm
Cl}_E(\mathcal{T})\not\subseteq\overline{\mathcal{T}}_{\rm
fin}\cup\overline{\mathcal{T}}_{\omega,1}$;

$(2)$ for some natural $n\geq 1$, Boolean algebras
$\mathcal{B}_n(T)$, $T\in\mathcal{T}$, have unbounded
cardinalities and, moreover, there is an infinite decreasing chain
$(u_k)_{k\in\omega}$ of labels for some formulas \
$\varphi_k(\bar{x})$ \ such that any finite sequence \
$(u_l,\ldots,u_0)$ with $u_l<\ldots<u_0$ is witnessed by
infinitely many theories in $\mathcal{T}$;

$(3)$ the same as in $(2)$ with $u_0=1$.}

\medskip
{\bf Corollary 2.2.} {\em A class
$\mathcal{T}\subseteq\overline{\mathcal{T}}_{\rm
fin}\cup\overline{\mathcal{T}}_{\omega,1}$ does not generate,
using the $E$-operator, theories, which are neither finitely
categorical and $\omega$-categorical, if and only if for the
Boolean algebras $\mathcal{B}_n(T)$, $T\in\mathcal{T}$, there are
no infinite decreasing chains $(u_k)_{k\in\omega}$ of labels for
some formulas $\varphi_k(\bar{x})$ such that any finite sequence
$(u_l,\ldots,u_0)$ with $u_l<\ldots<u_0$ is witnessed by
infinitely many theories in $\mathcal{T}$.}

\medskip
{\bf Remark 2.3.} Corollary 2.2 together with Proposition 1.1
allow to determine $E$-closed classes of finitely categorical and
$\omega$-categorical theories. Here, since finite sets of theories
are $E$-closed, it suffices to consider infinite sets.

\medskip
Considering a set $\mathcal{T}$ of theories with disjoint
languages, for the $E$-closeness it suffices to add theories of
the empty language describing cardinalities, in $\omega+1$, of
universes if these cardinalities meet infinitely many times in
$\mathcal{T}$.

In such a case we obtain relative closures \cite{rest} and have
the following assertions.

\medskip
{\bf Proposition 2.4.} {\em A class $\mathcal{T}$ of theories of
pairwise disjoint languages is $E$-closed if and only if the
following conditions hold:

{\rm (i)} for any $n\in\omega\setminus\{0\}$ whenever
$\mathcal{T}$ contains infinitely many theories with $n$-element
models then $\mathcal{T}$ contains the theory $T^0_n$ of the empty
language and with $n$-element models;

{\rm (ii)} if $\mathcal{T}$ contains theories with unbounded
finite cardinalities of models, or infinitely many theories with
infinite models, then $\mathcal{T}$ contains the theory
$T^0_\infty$ of the empty language and with infinite models.}

\medskip
{\bf Corollary 2.5.} {\em A class
$\mathcal{T}\subset\overline{\mathcal{T}}_{\rm
fin}\cup\overline{\mathcal{T}}_{\omega,1}$ of theories of pairwise
disjoint languages is $E$-closed if and only if the conditions
{\rm (i)} and {\rm (ii)} hold.}

\medskip
{\bf Corollary 2.6.} {\em A class
$\mathcal{T}\subset\overline{\mathcal{T}}_{\rm fin}$ of theories
of pairwise disjoint languages is $E$-closed if and only if the
condition {\rm (i)} holds and there is $N\in\omega$ such
$\mathcal{T}$ does not have $n$-categorical theories for $n>N$.}

\medskip
{\bf Corollary 2.7.} {\em A class
$\mathcal{T}\subset\overline{\mathcal{T}}_{\omega,1}$ of theories
of pairwise disjoint languages is $E$-closed if and only if the
condition {\rm (ii)} holds.}

\medskip
{\bf Corollary 2.8.} {\em For any class
$\mathcal{T}\subset\overline{\mathcal{T}}_{\rm
fin}\cup\overline{\mathcal{T}}_{\omega,1}$ of theories of pairwise
disjoint languages, ${\rm Cl}_E(\mathcal{T})$ is contained in the
class $\subset\overline{\mathcal{T}}_{\rm
fin}\cup\overline{\mathcal{T}}_{\omega,1}$, moreover,
$${\rm
Cl}_E(\mathcal{T})\subseteq\mathcal{T}\cup\{T^0_\lambda\mid\lambda\in(\omega\setminus\{0\})\cup\{\infty\}
\}.$$}

{\bf Remark 2.9.} Using relative closures \cite{rest} the
assertions 2.4--2.8 also hold if languages are disjoint modulo a
common sublanguage $\Sigma_0$ such that all restrictions of
$n$-categorical theories in
$\mathcal{T}\cap\overline{\mathcal{T}}_{\rm fin}$ to $\Sigma_0$
have isomorphic (finite) models $\mathcal{M}_n$ and all
restrictions of theories in
$\mathcal{T}\cap\overline{\mathcal{T}}_{\omega,1}$ have isomorphic
countable models $\mathcal{M}_\omega$. In such a case, the
theories $T^0_n$ should be replaced by ${\rm Th}(\mathcal{M}_n)$
and $T^0_\infty$~--- by ${\rm Th}(\mathcal{M}_\omega)$.

It is also permitted to have finitely many possibilities for each
$\mathcal{M}_n$ and for $\mathcal{M}_\omega$.

\medskip
The following example shows that (even with pairwise disjoint
languages) $\omega$-categorical theories $T$ with unbounded
$\rho_{n,T}(0,1)$ do not force theories outside the class of
$\omega$-categorical theories.

\medskip
{\bf Example 2.10.} Let $T_n$ be a theory of infinitely many
disjoint $n$-cubes with a graph relation $R^{(2)}_n$, $R_m\ne R_n$
for $m\ne n$. For the set $\mathcal{T}=\{T_n\mid n\in\omega\}$ we
have ${\rm Cl}_E(\mathcal{T})=\mathcal{T}\cup\{T^0_\infty\}$. All
theories in ${\rm Cl}_E(\mathcal{T})$ are $\omega$-categorical
whereas $\rho_{2,T_n}(0,1)=n+2$ that witnessed by formulas
describing distances $d(x,y)\in\omega\cup\{\infty\}$ between
elements.

Similarly, taking for each $n\in\omega$ exactly one $n$-cube with
a graph relation $R^{(2)}_n$, we get a set $\mathcal{T}$ of
theories such that ${\rm
Cl}_E(\mathcal{T})\subset\overline{\mathcal{T}}_{\rm
fin}\cup\overline{\mathcal{T}}_{\omega,1}$.

\medskip
{\bf Remark 2.11.} Assertions 2.1 -- 2.5 and 2.7 -- 2.9 hold for
the operators ${\rm Cl}^{d}_P$ and ${\rm Cl}^{d,r}_P$ replacing
$E$-closures by $P$-closures. As non-isolated types always produce
infinite structures, Corollary 2.6 holds only for ${\rm Cl}^{d}_P$
with finite sets $\mathcal{T}$ of theories.

\section{On approximations of theories with (in)finite models}

{\bf Definition} \cite{Rosen}. An infinite structure $\mathcal{M}$
is {\em pseudofinite} if every sentence true in $\mathcal{M}$ has
a finite model.

\medskip
{\bf Definition} (cf. \cite{strmin}). A consistent formula
$\varphi$ {\em forces} the infinity if $\varphi$ does not have
finite models.

\medskip
By the definition, an infinite structure $\mathcal{M}$ is
pseudofinite if and only if $\mathcal{M}$ does not satisfy
formulas forcing the infinity.

\medskip
We denote the class
$\overline{\mathcal{T}}\setminus\overline{\mathcal{T}}_{\rm fin}$
by $\overline{\mathcal{T}}_{\rm inf}$.

\medskip
{\bf Proposition 3.1.} {\em A theory
$T\in\overline{\mathcal{T}}_{\rm inf}$ belongs to some $E$-closure
of theories in $\overline{\mathcal{T}}_{\rm fin}$ if and only if
$T$ does not have formulas forcing the infinity.}

\medskip
{\bf\em Proof.} If a formula $\varphi$ forces the infinity then
$\mathcal{T}_\varphi\subset\overline{\mathcal{T}}_{\rm inf}$ for
any $\mathcal{T}\subseteq\overline{\mathcal{T}}$. Thus, having
such a formula $\varphi\in T$, $T$ can not be approximated by
theories in $\overline{\mathcal{T}}_{\rm fin}$ and so $T$ does not
belong to $E$-closures of families
$\mathcal{T}\subseteq\overline{\mathcal{T}}_{\rm fin}$.

Conversely, if any formula $\varphi\in T$ does not force the
infinity then, since $T\notin\mathcal{T}_{\rm fin}$,
$(\overline{\mathcal{T}}_{\rm fin})_\varphi$ is infinite using
unbounded finite cardinalities and we can choose infinitely many
theories in $(\overline{\mathcal{T}}_{\rm fin})_\varphi$, for each
$\varphi\in T$, forming a set
$\mathcal{T}_0\subset\overline{\mathcal{T}}_{\rm fin}$ such that
$T\in{\rm Cl}_E(\mathcal{T}_0)$.~$\Box$

\medskip
Note that, in view of Proposition 1.1, Proposition 3.1 is a
reformulation of Lemma 1 in \cite{Vaa}.

\medskip
{\bf Corollary 3.2.} {\em If a theory
$T\in\overline{\mathcal{T}}_{\rm inf}$ belongs to some $E$-closure
of theories in $\overline{\mathcal{T}}_{\rm fin}$ then $T$ is not
finitely axiomatizable.}

\medskip
{\bf\em Proof.} If $T$ is finitely axiomatizable by some formula
$\varphi$ then $|\mathcal{T}_\varphi|\leq 1$ for any
$\mathcal{T}\subseteq\overline{\mathcal{T}}$ and $\varphi$ forces
the infinity. Thus, in view of Proposition 3.1, $T$ can not be
approximated by theories in $\overline{\mathcal{T}}_{\rm fin}$,
i.~e., $T$ does not belong to $E$-closures of families
$\mathcal{T}_0\subset\overline{\mathcal{T}}_{\rm fin}$.~$\Box$

\medskip
In fact, in view of Theorem 1.2, the arguments for Corollary 3.3
show that ${\rm Cl}_E(\mathcal{T})$, for a family $\mathcal{T}$ of
finitely axiomatizable theories, has the least generating set
$\mathcal{T}$ and does not contain new finitely axiomatizable
theories.

\medskip
Note that Proposition 3.1 admits a reformulation for ${\rm
Cl}^d_P$ repeating the proof. At the same time theories in
$\overline{\mathcal{T}}_{\rm fin}$ can not be approximated by
theories in $\overline{\mathcal{T}}_{\rm inf}$ with respect to
${\rm Cl}_E$ (in view of Proposition 1.1) whereas each theory in
$\overline{\mathcal{T}}_{\rm fin}$ can be approximated by theories
in $\overline{\mathcal{T}}_{\rm inf}$ with respect to ${\rm
Cl}^d_P$:

\medskip
{\bf Proposition 3.3.} {\em For any theory
$T\in\overline{\mathcal{T}}_{\rm fin}$ there is a family
$\mathcal{T}_0\subset\overline{\mathcal{T}}_{\rm inf}$ such that
$T$ belongs to the $\Sigma(T)$-restriction of ${\rm
Cl}^d_P(\mathcal{T}_0)$.}

\medskip
{\bf\em Proof.} It suffices to form $\mathcal{T}_0$ by infinitely
many theories of structures $\mathcal{A}_i$, $i\in I$, with
infinitely many copies of models $\mathcal{M}\models T$ forming
$E_i$-classes for equivalence relations $E_i$, where $E_j$ is
either equality or complete for $j\ne i$. Considering disjoint
unary predicates $P_i$ for $\mathcal{A}_i$ we get the nonprincipal
$1$-type $p_\infty(x)$ isolated by the set $\{\neg P_i(x)\mid i\in
I\}$ which can be realized by the set $M$ with the structure
$\mathcal{M}$ witnessing that $T$ belongs to the restriction of
${\rm Cl}^d_P(\mathcal{T}_0)$ removing new relations $E_i$.~$\Box$

\medskip
{\bf Remark 3.4.} We have a similar effect removing all relations
$E_j$ in the structures $\mathcal{A}_i$ and obtaining isomorphic
structures $\mathcal{A}'_i$: by compactness the $P$-combination of
structures $\mathcal{A}'_i$ (where disjoint $\mathcal{A}'_i$ form
unary predicates $P_i$) has the theory with a model, whose
$p_\infty$-restriction forms a structure isomorphic to
$\mathcal{M}$. In this case we have ${\rm Cl}^{d,r}_P(\{{\rm
Th}(\mathcal{A}'_i)\})$.~$\Box$

\medskip
{\bf Remark 3.5.} As in the proof of Proposition 3.1 theories in
$\mathcal{T}_0$ can be chosen consistent modulo cardinalities of
their models we can add that $e$-${\rm Sp}(T)=1$ for the
$E$-combination $T$ of the theories in $\mathcal{T}_0$.

As the same time $e$-${\rm Sp}(T')$ is infinite for the
$P$-combination $T'$ of $\mathcal{A}_i$ in the proof of
Proposition 3.3, since $p_\infty(x)$ has infinitely many
possibilities for finite cardinalities of sets of realizations for
$p_\infty(x)$.~$\Box$

\section{$e$-spectra for finitely categorical \\ and $\omega$-ca\-te\-go\-ri\-cal theories}

We refine the notions of $e$-spectra $e$-${\rm Sp}(\mathcal{A}_E)$
and $e$-${\rm Sp}(T)$ for the theories $T={\rm Th}(\mathcal{A}_E)$
restricting the class of possible theories to a given class
$\mathcal{T}$ in the following way.

For a structure $\mathcal{A}_E$ the number of {\em
new}\index{Structure!new} structures with respect to the
structures $\mathcal{A}_i$, i.~e., of the structures $\mathcal{B}$
with ${\rm Th}(\mathcal{B})\in\mathcal{T}$, which are pairwise
elementary non-equivalent and elementary non-equivalent to the
structures $\mathcal{A}_i$, is called the {\em
$(e,\mathcal{T})$-spectrum}\index{$(e,\mathcal{T})$-spectrum} of
$\mathcal{A}_E$ and denoted by $(e,\mathcal{T})$-${\rm
Sp}(\mathcal{A}_E)$.\index{$(e,\mathcal{T})$-${\rm
Sp}(\mathcal{A}_E)$} The value ${\rm sup}\{(e,\mathcal{T})$-${\rm
Sp}(\mathcal{A}'))\mid\mathcal{A}'\equiv\mathcal{A}_E\}$ is called
the {\em
$(e,\mathcal{T})$-spectrum}\index{$(e,\mathcal{T})$-spectrum} of
the theory ${\rm Th}(\mathcal{A}_E)$ and denoted by
$(e,\mathcal{T})$-${\rm Sp}({\rm
Th}(\mathcal{A}_E))$.\index{$(e,\mathcal{T})$-${\rm Sp}({\rm
Th}(\mathcal{A}_E))$}

\medskip
The following properties are obvious.

\medskip
1. (Monotony) If $\mathcal{T}_1\subseteq\mathcal{T}_2$ then
$(e,\mathcal{T}_1)$-${\rm Sp}({\rm
Th}(\mathcal{A}_E))\leq(e,\mathcal{T}_2)$-${\rm Sp}({\rm
Th}(\mathcal{A}_E))$ for any structure $\mathcal{A}_E$.

\medskip
2. (Additivity) If the class $\overline{\mathcal{T}}$ of all
complete elementary theories of relational languages is the
disjoint union of subclasses $\overline{\mathcal{T}}_1$ and
$\overline{\mathcal{T}}_2$ then for any theory $T={\rm
Th}(\mathcal{A}_E)$, $$e\mbox{-}{\rm
Sp}(T)=(e,\overline{\mathcal{T}}_1)\mbox{-}{\rm
Sp}(T)+(e,\overline{\mathcal{T}}_2)\mbox{-}{\rm Sp}(T).$$

\medskip
We divide a class $\mathcal{T}$ of theories into two disjoint
subclasses $\mathcal{T}^{\rm fin}$\index{$\mathcal{T}^{\rm fin}$}
and $\mathcal{T}^{\rm inf}$\index{$\mathcal{T}^{\rm inf}$} having
finite and infinite non-empty language relations, respectively.
More precisely, for functions $f\mbox{\rm : }\omega\to\lambda_f$,
where $\lambda_f$ are cardinalities, we divide $\mathcal{T}$ into
subclasses $\mathcal{T}^f$ of theories $T$ such that $T$ has
$f(n)$ $n$-ary predicate symbols for each $n\in\omega$.

For the function $f$ we denote by ${\rm Supp}(f)$ its {\em
support}\index{Support}, i.e., the set $\{n\in\omega\mid
f(n)>0\}$.

\medskip
Clearly, the language of a theory $T\in\mathcal{T}^f$ is finite if
and only if $\rho_f\subset\omega$ and ${\rm Supp}(f)$ is finite.

\medskip
Illustrating $(e,\mathcal{T})$-spectra for the class $\mathcal{T}$
of all cubic theories and taking the class $\mathcal{T}^{\rm
fin}_0\subset\mathcal{T}$ of all theories of finite cubes we note
that for an $E$-combination $T$ of theories $T_i$ in
$\mathcal{T}^{\rm fin}_0$, $(e,\mathcal{T})\mbox{-}{\rm Sp}(T)$ is
positive if and only if there are infinitely many $T_i$. In such a
case, $(e,\mathcal{T})\mbox{-}{\rm Sp}(T)=1$ and new theory, which
does not belong to $\mathcal{T}^{\rm fin}_0$, is the theory of
$\omega$-cube.

\medskip
The class $\overline{\mathcal{T}}_{\rm fin}$ is represented as
disjoint union of subclasses $\overline{\mathcal{T}}_{\rm
fin,n}$\index{$\overline{\mathcal{T}}_{{\rm fin},n}$} of theories
having $n$-element models, $n\in\omega\setminus\{0\}$. For
$N\in\omega$, the class $\bigcup\limits_{n\leq
N}\overline{\mathcal{T}}_{{\rm fin},n}$ is denoted by
$\overline{\mathcal{T}}_{{\rm fin},\leq
N}$\index{$\overline{\mathcal{T}}_{{\rm fin}, \leq N}$}.

\medskip
{\bf Proposition 4.1.} {\em For any
$\mathcal{T}\subset\overline{\mathcal{T}}$, ${\rm
Cl}_E(\mathcal{T})\setminus\overline{\mathcal{T}}_{\rm
fin}\ne\varnothing$ if and only if for any natural $N$,
$\mathcal{T}\not\subset\overline{\mathcal{T}}_{{\rm fin}, \leq
N}$.}

\medskip
{\bf\em Proof.} If $\mathcal{T}$ contains a theory with infinite
models, the assertion is obvious. If
$\mathcal{T}\subset\overline{\mathcal{T}}_{\rm fin}$, then we
apply Compactness and Proposition 1.1.~$\Box$

\medskip
The following obvious proposition is also based on Proposition
1.1.

\medskip
{\bf Proposition 4.2.} {\em If
$\mathcal{T}\subset\overline{\mathcal{T}}_{{\rm fin},n}$ {\rm
(}respectively $\mathcal{T}\subset\overline{\mathcal{T}}_{{\rm
fin},\leq N}${\rm )} then ${\rm
Cl}_E(\mathcal{T})\subset\overline{\mathcal{T}}_{{\rm fin},n}$
{\rm (}${\rm Cl}_E(\mathcal{T})\subset\overline{\mathcal{T}}_{{\rm
fin},\leq N}${\rm )}. For any theory $T={\rm Th}(\mathcal{A}_E)$,
where all $E$-classes have theories in $\mathcal{T}$, $e$-${\rm
Sp}(T)=(e,\overline{\mathcal{T}}_{{\rm fin},n})$-${\rm Sp}({\rm
Th}(\mathcal{A}_E))$ {\rm (}$e$-${\rm
Sp}(T)=(e,\overline{\mathcal{T}}_{{\rm fin},\leq N})$-${\rm
Sp}({\rm Th}(\mathcal{A}_E))${\rm )}. If, additionally,
$\mathcal{T}$ is the set of theories in a finite language then
$\mathcal{T}$ is finite {\rm (}and so $E$-closed{\rm )}. In
particular, for any theory $T={\rm Th}(\mathcal{A}_E)$ in a finite
language, where all $E$-classes have theories in $\mathcal{T}$,
$e$-${\rm Sp}(T)=0$.}

\medskip
{\bf Remark 4.3.} In fact, the conclusions of Proposition 4.2
follow implying the following fact. If all theories in
$\mathcal{T}$ contain a formula $\varphi$ then all theories in
${\rm Cl}_E(\mathcal{T})$ contain $\varphi$. For (1) we take a
formula $\varphi$ ``saying'' that models have exactly $n$
elements, and for (2)~--- a formula $\varphi$ ``saying'' that
models have at most $N$ elements. If the language is finite there
are only finitely many possibilities for isomorphism types on
$n$-element sets and these possibilities are formula-definable.

\medskip
Similarly Proposition 4.2 we have

\medskip
{\bf Proposition 4.4.} {\em If
$\mathcal{T}\cap\overline{\mathcal{T}}_{\rm fin}=\varnothing$ then
${\rm Cl}_E(\mathcal{T})\cap\overline{\mathcal{T}}_{\rm
fin}=\varnothing$.}

\medskip
{\bf Definition} \cite{lut}. A theory $T$ in a predicate language
$\Sigma$ is called {\em language uniform},\index{Theory!language
uniform} or a {\em {\rm LU}-theory}\index{{\rm LU}-theory} if
for each arity $n$
any substitution on the set of non-empty $n$-ary predicates
(corresponding to the symbols in $\Sigma$) preserves $T$. The {\rm
LU}-theory $T$ is called {\em {\rm IILU}-theory}\index{{\rm
IILU}-theory} if it has non-empty predicates and as soon as there
is a non-empty $n$-ary predicate then there are infinitely many
non-empty $n$-ary predicates and there are infinitely many empty
$n$-ary predicates.

\medskip
Since for any finite cardinality $n$ there are {\rm IILU}-theories
with $n$-element models, repeating the proof of \cite[Proposition
12]{lut} and \cite[Proposition 13]{lut} we get

\medskip
{\bf Proposition 4.5.} {\em $(1)$ For any
$n\in\omega\setminus\{0\}$ and $\mu\leq\omega$ there is an
$E$-combination $T={\rm Th}(\mathcal{A}_E)$ of ${\rm
IILU}$-theories $T_i\in\overline{\mathcal{T}}_{{\rm fin},n}$ in a
language $\Sigma$ of the cardinality $\omega$ such that $T$ has an
$e$-least model and $e$-${\rm Sp}(T)=\mu$.

$(2)$ For any uncountable cardinality $\lambda$ there is an
$E$-combination $T={\rm Th}(\mathcal{A}_E)$ of ${\rm
IILU}$-theories $T_i\in\overline{\mathcal{T}}_{{\rm fin},n}$ in a
language $\Sigma$ of the cardinality $\lambda$ such that $T$ has
an $e$-least model and $e$-${\rm Sp}(T)=\lambda$.}

\medskip
{\bf Proposition 4.6.} {\em For any $n\in\omega\setminus\{0\}$ and
infinite cardinality $\lambda$ there is an $E$-combination $T={\rm
Th}(\mathcal{A}_E)$ of ${\rm IILU}$-theories
$T_i\in\overline{\mathcal{T}}_{{\rm fin},n}$ in a language
$\Sigma$ of cardinality $\lambda$ such that $T$ does not have
$e$-least models and $e$-${\rm Sp}(T)\geq{\rm
max}\{2^\omega,\lambda\}$.}

\medskip
{\bf Proposition 4.7.} {\em  For any $n\in\omega\setminus\{0\}$
and infinite cardinality $\lambda$ there is an $E$-combination
$T={\rm Th}(\mathcal{A}_E)$ of ${\rm LU}$-theories
$T_i\in\overline{\mathcal{T}}_{{\rm fin},n}$ in a language
$\Sigma$ of cardinality $\lambda$ such that $T$ does not have
$e$-least models and $e$-${\rm Sp}(T)=2^\lambda$.}

\medskip
{\bf\em Proof.} Let $\Sigma$ be a language consisting, for some
natural $m$, of $m$-ary predicate symbols $R_i$, $i<\lambda$. For
any $\Sigma'\subseteq\Sigma$ we take a structure
$\mathcal{A}_{\Sigma'}$ of the cardinality $n$ such that
$R_i=(A_{\Sigma'})^m$ for $R_i\in\Sigma'$, and $R_i=\varnothing$
for $R_i\in\Sigma\setminus\Sigma'$. Clearly, each structure
$\mathcal{A}_{\Sigma'}$ has a ${\rm LU}$-theory and
$\mathcal{A}_{\Sigma'}\not\equiv\mathcal{A}_{\Sigma''}$ for
$\Sigma'\ne\Sigma''$. For the $E$-combination $\mathcal{A}_E$ of
the structures $\mathcal{A}_{\Sigma'}$ we obtain the theory
$T={\rm Th}(\mathcal{A}_E)$ having a model of the cardinality
$\lambda$. At the same time $\mathcal{A}_E$ has $2^\lambda$
distinct theories of the $E$-classes $\mathcal{A}_{\Sigma'}$.
Thus, $e$-${\rm Sp}(T)=2^\lambda$. Finally we note that $T$ does
not have $e$-least models by Theorem 1.2 and arguments for
\cite[Proposition 9]{cl}.~$\Box$

\medskip
{\bf Remark 4.8.} Considering countable ${\rm LU}$-theories for
the assertions above we can assume that these theories belong to a
class $\mathcal{T}^f$, where $f\in\omega^\omega$ and ${\rm
Supp}(f)$ is infinite. Note also that Propositions 4.5--4.7 hold
replacing the classes $\overline{\mathcal{T}}_{{\rm fin},n}$ by
$\overline{\mathcal{T}}_{\omega,1}$.

\medskip
Replacing $E$-classes by unary predicates $P_i$ (not necessary
disjoint) being universes for structures $\mathcal{A}_i$ and
restricting models of ${\rm Th}(\mathcal{A}_P)$ to the set of
realizations of $p_\infty(x)$ we get the {\em
$(e,\mathcal{T})$-spectrum}\index{$(e,\mathcal{T})$-spectrum}
$(e,\mathcal{T})$-${\rm Sp}({\rm
Th}(\mathcal{A}_P))$\index{$(e,\mathcal{T})$-${\rm Sp}({\rm
Th}(\mathcal{A}_P))$}, i.~e., the number of pairwise elementary
non-equivalent restrictions $\mathcal{N}$ of
$\mathcal{M}\models{\rm Th}(\mathcal{A}_P)$ to $p_\infty(x)$ such
that ${\rm Th}(\mathcal{N})\in\mathcal{T}$.

\medskip
{\bf Proposition 4.9.} {\em If the structures $\mathcal{A}_i$ have
pairwise disjoint languages with disjoint predicates $P_i$ then
for any natural $n\geq 1$, $(e,\overline{\mathcal{T}}_{{\rm
fin},n})$-${\rm Sp}({\rm Th}(\mathcal{A}_P))\leq 1$, and
$(e,\overline{\mathcal{T}}\setminus\overline{\mathcal{T}}_{{\rm
fin}})$-${\rm Sp}({\rm Th}(\mathcal{A}_P))\leq 1$.}

\medskip
{\bf\em Proof.} Clearly, if the structures $\mathcal{A}_i$ have
pairwise disjoint languages with disjoint predicates $P_i$ then
structures for $p_\infty(x)$ do not contain realizations of
language predicates, i.~e., have theories $T^0_\lambda$. Now
$(e,\overline{\mathcal{T}}_{{\rm fin},n})$-${\rm Sp}({\rm
Th}(\mathcal{A}_P))\leq 1$ and $(e,\overline{\mathcal{T}}_{{\rm
fin},n})$-${\rm Sp}({\rm Th}(\mathcal{A}_P))= 1$ if and only if
there are infinitely many indexes $i$ and ${\rm
Th}(\mathcal{A}_i)\ne T^0_n$ for any $i$. Similarly,
$(e,\overline{\mathcal{T}}\setminus\overline{\mathcal{T}}_{{\rm
fin}})$-${\rm Sp}({\rm Th}(\mathcal{A}_P))\leq 1$ and
$(e,\overline{\mathcal{T}}\setminus\overline{\mathcal{T}}_{{\rm
fin}})$-${\rm Sp}({\rm Th}(\mathcal{A}_P))= 1$ if and only if
there are infinitely many indexes $i$ and ${\rm
Th}(\mathcal{A}_i)\ne T^0_\infty$ for any $i$.~$\Box$

\medskip
Clearly, approximating structures without non-trivial predicates
and applying the proof of Proposition 4.9 we get a family of
$P$-combinations with $(e,\overline{\mathcal{T}}_{{\rm
fin},n})$-${\rm Sp}({\rm Th}(\mathcal{A}_P))=1$, for
$n\in\omega\setminus\{0\}$, and
$(e,\overline{\mathcal{T}}\setminus\overline{\mathcal{T}}_{{\rm
fin}})$-${\rm Sp}({\rm Th}(\mathcal{A}_P))=1$.

\medskip
Comparing approximations in Section 3 and proofs for
\cite[Propositions 4.12, 4.13]{cs} we get

\medskip
{\bf Proposition 4.10.} {\em For any infinite cardinality
$\lambda$ there is a theory $T={\rm Th}(\mathcal{A}_P)$ being a
$P$-combination of theories in $\overline{\mathcal{T}}_{{\rm
fin}}$ and of a language $\Sigma$ such that $|\Sigma|=\lambda$ and
$e$-${\rm Sp}(T)=2^\lambda$.}

\section{Almost language uniform theories}

{\bf Definition.} A theory $T$ in a predicate language $\Sigma$ is
called {\em almost language uniform},\index{Theory!almost language
uniform} or a {\em {\rm ALU}-theory}\index{{\rm ALU}-theory} if
for each arity $n$ with $n$-ary predicates for $\Sigma$ there is a
partition for all $n$-ary predicates, corresponding to the symbols
in $\Sigma$, with finitely many classes $K$ such that
any substitution preserving these classes preserves $T$, too. The
{\rm ALU}-theory $T$ is called {\em {\rm IIALU}-theory}\index{{\rm
IIALU}-theory} if it has non-empty predicates and as soon as there
is a non-empty $n$-ary predicate in a class $K$ then there are
infinitely many non-empty $n$-ary predicates in $K$ and there are
infinitely many empty $n$-ary predicates.

\medskip
By the definition any ${\rm LU}$-theory is an {\rm ALU}-theory and
any ${\rm IILU}$-theory is an {\rm IIALU}-theory as well.

Since any finite structure can have only finitely many distinct
predicates for each arity $n$ we get the following

\medskip
{\bf Proposition 5.1.} {\em Any theory
$T\in\overline{\mathcal{T}}_{\rm fin}$ is an {\rm ALU}-theory.}

\medskip
Replacing ${\rm LU}$- and ${\rm IILU}$- by {\rm ALU}- and {\rm
IIALU}- and the proofs in Propositions 4.5--4.7 we get analogs for
these assertions attracting expansions of arbitrary theories in
$\overline{\mathcal{T}}_{{\rm fin},n}$. Thus any theory in
$\overline{\mathcal{T}}_{{\rm fin},n}$ can be used obtaining
described $e$-spectra.

\section{Families of cardinalities for models of theories in closures}

Let $\mathcal{T}$ be a nonempty family of theories in
$\overline{\mathcal{T}}$. We denote by $c_E(\mathcal{T})$
(respectively, $c_P(\mathcal{T})$, $c^d_P(\mathcal{T})$,
$c^{d,r}_P(\mathcal{T})$) the set of finite cardinalities for
models of theories in ${\rm Cl}_E(\mathcal{T})$ (${\rm
Cl}_P(\mathcal{T})$, ${\rm Cl}^d_P(\mathcal{T})$, ${\rm
Cl}^{d,r}_P(\mathcal{T})$) and by $\bar{c}_E(\mathcal{T})$
(respectively, $\bar{c}_P(\mathcal{T})$,
$\bar{c}^d_P(\mathcal{T})$, $\bar{c}^{d,r}_P(\mathcal{T})$) the
set of finite cardinalities for models of theories in ${\rm
Cl}_E(\mathcal{T})$ (${\rm Cl}_P(\mathcal{T})$, ${\rm
Cl}^d_P(\mathcal{T})$, ${\rm Cl}^{d,r}_P(\mathcal{T})$) which are
not cardinalities for models of theories in $\mathcal{T}$.
Additionally, for ${\rm Cl}_P(\mathcal{T})$, ${\rm
Cl}^d_P(\mathcal{T})$ and ${\rm Cl}^{d,r}_P(\mathcal{T})$ we
denote by $\hat{c}_P(\mathcal{T})$, $\hat{c}^d_P(\mathcal{T})$,
$\hat{c}^{d,r}_P(\mathcal{T})$, respectively, the set of finite
cardinalities for models of theories being restrictions for
corresponding $P$-combinations to sets of realizations of types
$p_\infty(x)$.

\medskip
{\bf Remark 6.1.} Since $E$-closures preserve finite cardinalities
for models of theories in families in $\mathcal{T}$, i.e.,
$c_E(\mathcal{T})$ consists of these cardinalities for
$\mathcal{T}$, then $\bar{c}_E(\mathcal{T})\equiv\varnothing$.
Thus we can use the notation $c_E(\mathcal{T})$ for the set of
finite cardinalities for models of theories in $\mathcal{T}$, or,
equivalently, for models of theories in ${\rm Cl}_E(\mathcal{T})$.

\medskip
{\bf Remark 6.2.} If $\mathcal{T}$ is finite, or corresponding
$p_\infty(x)$ is consistent and there are no models with finitely
many realizations for $p_\infty(x)$, then
$c_P(\mathcal{T})=c^d_P(\mathcal{T})=c_E(\mathcal{T})$ and
$\bar{c}_P(\mathcal{T})=\bar{c}^d_P(\mathcal{T})=\hat{c}_P(\mathcal{T})=\hat{c}^d_P(\mathcal{T})=\varnothing$.

\medskip
Examples of families of theories in the empty language $\Sigma_0$
witness that the cardinalities for sets of realizations of
$p_\infty(x)$ can vary arbitrarily and for finite $\mathcal{T}$ we
have
$c^{d,r}_P(\mathcal{T})=\hat{c}^{d,r}_P(\mathcal{T})=\mathbb{Z}^+$
and $\bar{c}^{d,r}_P(\mathcal{T})=\mathbb{Z}^+\setminus
c_E(\mathcal{T})$.

Having an infinite family $\mathcal{T}$ in the language
$\Sigma_0$, similarly we get
$c_P(\mathcal{T})=c^d_P(\mathcal{T})=c^{d,r}_P(\mathcal{T})=\hat{c}_P(\mathcal{T})=\hat{c}^d_P(\mathcal{T})=\hat{c}^{d,r}_P(\mathcal{T})=\mathbb{Z}^+$
and
$\bar{c}_P(\mathcal{T})=\bar{c}^d_P(\mathcal{T})=\bar{c}^{d,r}_P(\mathcal{T})=\mathbb{Z}^+\setminus
c_E(\mathcal{T})$. The latter formula shows that
$\bar{c}_P(\mathcal{T})$, $\bar{c}^d_P(\mathcal{T})$, and
$\bar{c}^{d,r}_P(\mathcal{T})$ can be arbitrary subsets of
$\mathbb{Z}^+$ with infinite complements. Thus we have the
following

\medskip
{\bf Proposition 6.3.} {\em For any infinite set
$Y\subseteq\mathbb{Z}^+$ there is a family $\mathcal{T}$ such that
$\bar{c}_P(\mathcal{T})=\bar{c}^d_P(\mathcal{T})=
\bar{c}^{d,r}_P(\mathcal{T})=\mathbb{Z}^+\setminus Y$.}

\medskip
{\bf Example 6.4.} If the language $\Sigma$ consists of the symbol
$E_k$ of the equivalence relation whose each class has
$k\in\omega$ elements then $p_\infty(x)$ can form an arbitrary
structure with $k$-element equivalence classes and for a finite
family $\mathcal{T}_k$ we have
$c^{d,r}_P(\mathcal{T}_k)=\hat{c}^{d,r}_P(\mathcal{T}_k)=k\mathbb{Z}^+$
and $\bar{c}^{d,r}_P(\mathcal{T}_k)=k\mathbb{Z}^+\setminus
c_E(\mathcal{T}_k)$. If the family $\mathcal{T}_k$ is infinite
then, similarly,
$c_P(\mathcal{T}_k)=c^d_P(\mathcal{T}_k)=c^{d,r}_P(\mathcal{T}_k)=\hat{c}_P(\mathcal{T}_k)=\hat{c}^d_P(\mathcal{T}_k)=\hat{c}^{d,r}_P(\mathcal{T}_k)=k\mathbb{Z}^+$
and
$\bar{c}_P(\mathcal{T}_k)=\bar{c}^d_P(\mathcal{T}_k)=\bar{c}^{d,r}_P(\mathcal{T}_k)=k\mathbb{Z}^+\setminus
c_E(\mathcal{T})$.

More generally, collecting the families of theories with distinct
$E_k$, $k\in K$, $K\subset\omega$, we obtain nonempty values for
$c_P$, $c^d_P$, $c^{d,r}_P$, $\hat{c}_P$, $\hat{c}^d_P$,
$\hat{c}^{d,r}_P$ as unions $\bigcup\limits_{k\in
K}k\mathbb{Z}^+$.

\medskip
Now we have to show that all possible nonempty values for
$\hat{c}^d_P$ and $\hat{c}^{d,r}_P$ are exhausted by the sums
$\biguplus\limits_{k\in K}k\mathbb{Z}^+$ (unions with finite sums
for numbers in $k\mathbb{Z}^+$) whereas values for $\hat{c}_P$ may
differ.

\medskip
{\bf Theorem 6.5.} {\em For any nonempty family $\mathcal{T}$
there is $K\subset\omega$ such that
$\hat{c}^{d,r}_P(\mathcal{T})=\biguplus\limits_{k\in
K}k\mathbb{Z}^+$.}

\medskip
{\bf\em Proof.} Recall that for $P$-combinations with respect to
${\rm Cl}^{d,r}_P$ there are no links between disjoint predicates
$P_i$ with structures $\mathcal{A}_i$ being models of theories in
$\mathcal{T}$. Therefore if $p_\infty(x)$ can produce finite
structures then structures $\mathcal{A}_i$ with $1$-types
approximating $p_\infty(x)$, define (partial) definable
equivalence relations with bounded finite classes $E(a)$ and
without definable extensions, for the approximations and for
$p_\infty(x)$. So there are no links between the classes $E(a)$
and having $k$ elements in $E(a)$ we produce, by compactness, a
series of $1,2,\ldots,n,\ldots$ $E$-classes for $p_\infty(x)$
since $p_\infty(x)$ is not isolated. Thus we get a series
$k\mathbb{Z}^+$ for $\hat{c}^{d,r}_P(\mathcal{T})$. Varying finite
cardinalities for the classes $E(a)$ we obtain the required
formula $\hat{c}^{d,r}_P(\mathcal{T})=\biguplus\limits_{k\in
K}k\mathbb{Z}^+$ for some set $K\subset\omega$ witnessing these
cardinalities. If $p_\infty(x)$ can produce finite structures then
we set $K\rightleftharpoons\varnothing$.~$\Box$

\medskip
{\bf Theorem 6.6.} {\em For any infinite family $\mathcal{T}$
there is $K\subset\omega$ such that
$\hat{c}^{d}_P(\mathcal{T})=\biguplus\limits_{k\in
K}k\mathbb{Z}^+$.}

\medskip
{\bf\em Proof} repeats the proof of Theorem 6.5 using structures
$\mathcal{A}_i$ which pairwise are not elementary
equivalent.~$\Box$

\medskip
{\bf Remark 6.7.} 1. In Theorems 6.5 and 6.6, if we have minimal
$K$ with $|K|>1$ then the type $p_\infty(x)$ is not complete.
Indeed, taking, for sets of realizations of $p_\infty(x)$, maximal
definable equivalence relations $E_1$ and $E_2$ for $k_1\ne k_2\in
K$ we can not move, by automorphisms, elements of $E_1$-classes to
elements of $E_2$-classes.

2. Clearly, having $E_1$-classes and $E_2$-classes of same
cardinalities with non-isomorphic structures we again can not
connect elements of these classes by automorphisms. Thus, $|K|=1$
is a necessary but not sufficient condition for the completeness
of $p_\infty(x)$.

3. The least cardinality $|K|$, with positive $\hat{c}^{d,r}_P$ or
$\hat{c}^{d}_P$, gives a lower bound for independent equivalence
relations with respect to their realizability/omitting for
restrictions of models to sets of realizations of
$p_\infty(x)$.~$\Box$

\medskip
{\bf Remark 6.8.} Finite structures $\mathcal{A}_\infty$ for
maximal definable equivalence relations for $p_\infty(x)$ with
respect to ${\rm Cl}^d_P$ and to ${\rm Cl}^{d,r}_P$ can be
isomorphic if and only if they are represented in some
$\mathcal{A}_i$ for ${\rm Cl}^{d,r}_P$ and infinitely many times
for ${\rm Cl}^d_P$, or approximated both for ${\rm Cl}^{d,r}_P$
and for ${\rm Cl}^d_P$. Hence, for any infinite family
$\mathcal{T}$,
$\hat{c}^{d}_P(\mathcal{T})=\hat{c}^{d,r}_P(\mathcal{T})$ if and
only if each $n$-element class for maximal definable equivalence
relations for $p_\infty(x)$ with respect to ${\rm Cl}^{d,r}_P$ has
$n$-element classes for correspondent definable equivalence
relations in infinitely many pairwise elementary non-equivalent
structures $\mathcal{A}_i$, with respect to ${\rm
Cl}^{d}_P$.~$\Box$

\medskip
{\bf Definition} \cite{Pi83}. Let ${\cal M}$ be a model of a
theory~$T$, $\bar{a}$ and $\bar{b}$ tuples in ${\cal M}$, $A$ a
subset of $M$. The tuple $\bar{a}$
\emph{semi-isolates}\index{Tuple!semi-isolates a tuple} the tuple
$\bar{b}$ over the set $A$ if there exists a formula
$\varphi(\bar{a},\bar{y})\in{\rm tp}(\bar{b}/A\bar{a})$ for which
$\varphi(\bar{a},\bar{y})\vdash{\rm tp}(\bar{b}/A)$ holds. In this
case we say that the formula $\varphi(\bar{x},\bar{y})$ (with
parameters in $A$) \emph{witnesses that $\bar{b}$ is semi-isolated
over $\bar{a}$ with respect to
$A$}.\index{Formula!witnessing!semi-isolation}

If $p\in S(T)$ and ${\cal M}\models T$ then ${\rm SI}^{\cal
M}_p$\index{${\rm SI}^{\cal M}_p$} denotes the relation of
semi-isolation (over $\varnothing$) on the set of all realizations
of $p$:
$$
{\rm SI}^{\cal M}_p\rightleftharpoons\{(\bar{a},\bar{b})\mid{\cal
M}\models p(\bar{a})\wedge p(\bar{b}) \mbox{ and }\bar{a}\mbox{
semi-isolates }\bar{b}\}.
$$

\medskip
The following definition generalizes the previous one for a family
of $1$-types, in particular, for incomplete $p_\infty(x)$.

\medskip
{\bf  Definition} \cite{Su131}. Let $T$ be a complete theory,
$\mathcal{M}\models T$. We consider {\em closed}\index{Set!closed}
nonempty sets (under the natural topology) sets ${\bf
p}(x)\subseteq S^1(\varnothing)$, i.~e., sets ${\bf p}(x)$ such
that ${\bf p}(x)=\bigcap\limits_{i\in I}[\varphi_{{\bf p},i}(x)]$,
where $[\varphi_{{\bf p},i}(x)]\rightleftharpoons\{p(x)\in
S^1(\varnothing)\mid \varphi_{{\bf p},i}(x)\in
p(x)\}$\index{$[\varphi_{{\bf p},i}(x)]$} for some formulas
$\varphi_{{\bf p},i}(x)$ of $T$.

For closed sets ${\bf p}(x),{\bf q}(y)\subseteq S(\varnothing)$ of
types, realized in $\mathcal{M}$, we consider {\em $({\bf p},{\bf
q})$-preserving}\index{Formula!$({\bf p},{\bf q})$-preserving}
{\em $({\bf p},{\bf q})$-semi-isolating}\index{Formula!$({\bf
p},{\bf q})$-semi-isolating}, {\em $({\bf p}\rightarrow {\bf
q})$-}\index{$({\bf p}\rightarrow{\bf q})$-formulas}, or {\em
$({\bf q}\leftarrow {\bf p})$-formulas}\index{$({\bf q}\leftarrow
{\bf p})$-formula} $\varphi(x,y)$ of $T$, i.~e., formulas for
which if $a\in M$ realizes a type in ${\bf p}(x)$ then every
solution of $\varphi(a,y)$ realizes a type in ${\bf q}(y)$.

If ${\bf p}(x)={\bf q}(y)$ then $({\bf p},{\bf q})$-preserving
formulas are called {\em ${\bf p}$-preserving}\index{Formula!${\bf
p}$-preserving} or {\em ${\bf
p}$-semi-isolating}\index{Formula!${\bf p}$-semi-isolating} and we
define, similarly to $ {\rm SI}^{\cal M}_p$, the {\em generalized}
relation ${\rm SI}^{\cal M}_{\bf p}$\index{${\rm SI}^{\cal M}_{\bf
p}$} of semi-isolation for the set of realizations of types in
${\bf p}(x)$:

$$
{\rm SI}^{\cal M}_{\bf p}\rightleftharpoons\{(a,b)\mid{\cal
M}\models p(a)\wedge p'(b)\wedge\varphi(a,b)$$
$$
\mbox{ for }p,p'\in{\bf p}\mbox{ and a }{\bf p}\mbox{-preserving
formula }\varphi(x,y)\}.
$$

If $(a,b)\in{\rm SI}^{\cal M}_{\bf p}$ we say that $a$ {\em
semi-isolates} $b$ {\em with respect to} ${\bf p}$.

\medskip
Thus, $a$ semi-isolates $b$ (in sense of \cite{Pi83}) if and only
if $a$ semi-isolates $b$ with respect to $\{{\rm tp}(a),{\rm
tp}(b)\}$.

\medskip
{\bf Remark 6.9.} Since there are no links between structures
$\mathcal{A}_i$ with respect to ${\rm Cl}^{d}_P$ and ${\rm
Cl}^{d,r}_P$, the set ${\bf p}_\infty$ of all completions $q(x)$
of $p_\infty(x)$ has symmetric ${\rm SI}^{\cal M}_{{\bf
p}_\infty}$. Thus, the relations ${\rm SI}^{\cal M}_{{\bf
p}_\infty}$ form equivalence relations. Positive values for
$\hat{c}^{d}_P$ and $\hat{c}^{d,r}_P$ imply that these equivalence
relations have finite classes. Cardinalities of these classes
define formulas in Theorems 6.5 and 6.6.~$\Box$

\medskip
Now we consider the general case, with the operator ${\rm Cl}_P$.
One can hardly expect productive descriptions considering
arbitrary links of structures with respect to arbitrary links of
predicates $P_i$, in contrast to the disjoint predicates when,
obviously, there are no links between the structures. So we will
fix a $P$-combination $\mathcal{A}_P$ (and its theory $T={\rm
Th}(\mathcal{A}_P)$) and consider the set $\hat{c}_P(T)$ of values
of finite cardinalities for $p_\infty(x)$ with respect to given
$P$-combination $T$, instead of the set $\hat{c}_P(\mathcal{T})$
of values for all finite values for all possible $P$-combinations.
In other words we argue to describe sets of finite cardinalities
for sets of realizations of a nonprincipal, not necessary
complete, $1$-type $p_\infty(x)$.

We note the following obvious observations.

\medskip
{\bf Remark 6.10.} 1. If any $n\in\omega$ realizations of a type
$p_\infty(x)$ force infinitely many realizations of $p_\infty(x)$
then it is true for any $m>n$.

2. If $a$ and $b$ are realizations of a type $p_\infty(x)$ and $a$
does not semi-isolate $b$ with respect to ${\bf p}_\infty$ then
there are no formulas $\varphi(x,y)$ with $\models\varphi(a,b)$
and forcing finitely or infinitely many realizations for the type
$q={\rm tp}(b/a)$, i.~e., the set of realizations of $q$ can be
empty and infinite, depending on a model.

The first observation shows that having $n$ which forces infinity,
we get $\hat{c}_P(T)\subset n$. The second one implies that
realizations of ${\bf p}_\infty$, which are not connected by the
relation of semi-isolation, contribute to $\hat{c}_P(T)$
independently on the binary level. Moreover, these contributions
by realizations $a$ and $b$ can generate distinct series, as in
Theorems 6.5 and 6.5, only if ${\rm tp}(a)\ne{\rm tp}(b)$.

\medskip
The following example shows that there is a theory $T$ with
$\hat{c}_P(T)=\{1\}$ clarifying that contributions above on the
binary level deny by the ternary level.

\medskip
{\bf Example 6.11.} Consider a coloring ${\rm Col}\mbox{\rm :
}M\to\omega\cup\{\infty\}$ of an infinite set $M$ such that each
color $\lambda\in\omega\cup\{\infty\}$ has infinitely many
elements in $M$, i.~e., each ${\rm Col}_n=\{a\in M\mid {\rm
Col}(a)=n\}$ is infinite as well as there are infinitely many
elements of the infinite color. We put
$P_i=M\setminus\bigcup\limits_{j<i}{\rm Col}_j$ and
$p_\infty(x)=\{\neg P_n(x)\mid n\in\omega\}$. Now we define, using
a generic construction with free amalgams \cite{SuCCMCT, Su073}, a
ternary relation $R$ such that for the definable relation
$Q(x,y)\equiv \exists z R(x,y,z)$ we have the following
properties:

1) the $Q$-structure has unique $1$-type and, moreover, its
automorphism group is transitive;

2) $R(x,y,z)\equiv Q(x,z)\wedge Q(y,z)$;

3) ${\rm Col}$ is an inessential coloring which is not neither
$Q$-ordered nor $Q^{-1}$-ordered \cite{SuCCMCT, Su031}, moreover,
for any element $a\in M$ the sets of solutions for $Q(a,y)$ and
$Q(x,a)$ have infinitely many elements of each color;

4) for any $a\ne b\in M$ the set of solutions for $R(a,b,z)$ is
infinite for each color $n\geq{\rm min}\{{\rm Col}(a),{\rm
Col}(b)\}$ and does not have have elements of colors $<{\rm
min}\{{\rm Col}(a),{\rm Col}(b)\}$, hence,  $R(a,b,z)\vdash
p_\infty(z)$ if $\models p_\infty(a)$ and $\models p_\infty(b)$.

Taking the generic structure $\mathcal{M}$ in the language
$\langle P^{(1)}_n, Q^{(2)},R^{(3)}\rangle_{n\in\omega}$ and its
theory $T={\rm Th}(\mathcal{M})$, being a $P$-combination, we have
$\hat{c}_P(T)=\{1\}$ since the nonisolated type $p_\infty(x)$ can
have, in a model of $T$, $0$, $1$, or infinitely many
realizations: one realization of $p_\infty(x)$ does not force new
ones and two distinct realizations $a$ and $b$ of $p_\infty(x)$
force infinitely many ones by the formula $R(a,b,z)$.

\medskip
{\bf Example 6.12.} We modify Example 6.11 replacing elements by
$E_k$-classes, where each class contains $k$ elements, and repeat
the generic construction satisfying the following conditions:

2) if $aE_k a'$ then ${\rm Col}(a)={\rm Col}(a')$;

2) if $(a,b)\in Q$, $aE_k a'$, $bE_kb'$, then $(a',b')\in Q$.

The theory $T_k$ of resulting generic structure $\mathcal{M}_k$
satisfies $\hat{c}_P(T_k)=\{k\}$ since each realization $a$ of
$p_\infty(x)$ forces $k$ realizations of $p_\infty(x)$ consisting
of $E_k(a)$ and any two realizations of $p_\infty(x)$ belonging to
distinct $E_k$-classes forces infinitely many $E_k$-classes with
elements satisfying $p_\infty(x)$.

\medskip
Combining structures $\mathcal{M}_k$ with distinct $k$ we obtain a
generic structure whose theory $T$ satisfies $\hat{c}_P(T)=K$ for
a given set $K\subseteq\mathbb Z^+$. Here sets of realizations of
$p_\infty(x)$ are divided into $E_k$-classes for $k\in K$.

Thus we have the following theorem asserting that values
$\hat{c}_P(T)$ can be arbitrary.

\medskip
{\bf Theorem 6.13.} {\em For any set $K\subseteq\mathbb Z^+$ there
is a $P$-combination $T$ such that $\hat{c}_P(T)=K$.}

\medskip
Now we argue to modify the generic construction and Theorem 6.13
using transitive arrangements of algebraic systems similar to
\cite{SuGP, Su1999}, and obtaining a similar result with complete
$p_\infty(x)$ describing possibilities for $\hat{c}_P(T)$.

For this aim we fix a nonempty set $K\subseteq\mathbb Z^+$
claiming for $\hat{c}_P(T)=K$ with some $P$-combination $T$. Note
that if $1\notin K$ then either any realization $a$ of
$p_\infty(x)$ forces infinitely many realizations or $a$ belongs
to the maximal finite definable $E$-class with some $k_0>1$
elements. At first case, by completeness of $p_\infty(x)$, any
finite set of realizations of $p_\infty(x)$ forces that infinity
and therefore $K=\varnothing$ contradicting the condition
$K\ne\varnothing$. At second case, again by completeness of
$p_\infty(x)$, we have $K\subseteq k_0\mathbb Z^+$. Replacing
elements by their $E$-classes we reduce the problem of
construction of $T$ with $\hat{c}_P(T)=K$ to the case $1\in K$.

Example 6.11 witnesses the possibility for $\hat{c}_P(T)=\{1\}$.
So below we assume that $1\in K$ and $|K|\geq 2$.

Now for each $k\in K\setminus\{1\}$ we introduce a ternary
relation $R_k$ defining a free (acyclic) precise pseudoplane
$\mathcal{P}_k$ \cite{SuGP} with infinitely many lines containing
any fixed point and exactly $k$ points belonging to any fixed line
such that $\mathcal{P}_k$ has infinitely many connected
components. Then we combine these free pseudoplanes
$\mathcal{P}_k$ allowing that each point belongs to each
pseudoplane $\mathcal{P}_k$ and the union of sets of lines does
not form cycles. We embed copies of that combination $\mathcal{P}$
of the pseudoplanes into unary predicates ${\rm Col}_n$ as well as
to the structure of $p_\infty(x)$.

Modifying Example 6.11 we introduce a binary predicate $Q$ such
that;

1) if $(a,b)\in Q$ and $(a,c)\in Q$ then $a$, $b$, $c$ belong to
pairwise distinct connected components of $\mathcal{P}$, the same
is satisfied for $Q^{-1}$ (as in Example described in
\cite[Section 1.3]{SuCCMCT} and in \cite{Su92});

2) elements $a_1,\ldots,a_m$, $m>1$, realizing $p_\infty(x)$ and
belonging to a common line $l$ force all elements of $l$ and do
not force elements outside $l$;

3) if $a$ and $b$ are realizations of $p_\infty(x)$ which do not
have a common line then $a$ and $b$ force infinitely many
realizations of $p_\infty(x)$ by the formula $Q(a,y)\wedge
Q(b,y)$.

The resulted generic structure $\mathcal{M}$ of the language
$\langle{\rm Col}_n,Q,R_k\rangle_{n\in\omega,k\in
K\setminus\{1\}}$ and its theory $T$ satisfy the following
properties:

i) any realization of $p_\infty(x)$ does not force new
realizations of $p_\infty(x)$ witnessing $1\in\hat{c}_P(T)$;

ii) any at least two distinct realizations of $p_\infty(x)$ in a
line $l$ belonging to $\mathcal{P}_k$ force exactly the set $l$
witnessing $k\in\hat{c}_P(T)$ for $k\in K$;

iii) any two distinct realizations of $p_\infty(x)$ which do not
have common lines force infinitely many realizations of
$p_\infty(x)$ witnessing $k'\notin\hat{c}_P(T)$ for $k'\notin K$.

Thus we get $\hat{c}_P(T)=K$.

\medskip
Collecting the arguments above we have the following

\medskip
{\bf Theorem 6.14.} {\em {\rm (1)} If $T$ is a $P$-combination
with a type $p_\infty(x)$ isolating a complete $1$-type then
$\hat{c}_P(T)$ is either empty or contains $k_0$ such that
$\hat{c}_P(T)\subseteq k_0\mathbb Z^+$.

{\rm (2)} For any set $K\subseteq k_0\mathbb Z^+$, being empty or
containing $k_0$, there is a $P$-combination $T$ with a type
$p_\infty(x)$ isolating a complete $1$-type such that
$\hat{c}_P(T)=K$.}


\end{document}